\documentclass{amsart}
\usepackage{graphicx}
\newtheorem{thm}{Theorem}[section]
\newtheorem{cor}[thm]{Corollary}
\newtheorem{lem}[thm]{Lemma}

\theoremstyle{definition}

\theoremstyle{remark}

\numberwithin{equation}{section}
\theoremstyle{example}

\begin{document}

 \title{ cr-hypersurfaces of a conformal Kenmotsu space form satisfying certain shape operator conditions}%
\author{R.Abdi and E.Abedi }%
\address{
Department of Mathematics\\
Azarbaijan shahid Madani University,Tabriz 53751 71379, Iran}%
\email{rabdi@azaruniv.edu}
\email{esabedi@azaruniv.edu}%



\begin{abstract}
In this paper, Conformal Kenmotsu manifolds are introduced which are not  Kenmotsu. We consider CR-hypersurfaces of a conformal Kenmotsu space form whose shape operator is parallel, scalar, recurrent or Lie $ \xi $-parallel,  it  is proved that if the Lee vector field of a conformal Kenmotsu space form is tangent  and normal to these type CR-hypersurfaces then the CR-hypersurfaces are totally geodesic and totally umbilic, respectively.
\end{abstract}
\maketitle
\section{Introduction}
Let $(M^{2n},J,g)$ be a Hermitian manifold of complex dimension n, where $J$ denotes its complex structure and $g$ its Hermitian metric. Then $(M^{2n},J,g)$ is a locally conformal
 K\"{a}ehler manifold if there is  an open cover $\{U_i\}_{i \in I}$ of $M^{2n}$ and a family $\{f_i\}_{i \in I}$ of $C^{\infty}$ functions $f_i: U_i \longrightarrow \mathbb{R}$ so that each local metric $g_{i} = exp(-f_i)g_{\mid U_i}$ is K\"{a}ehlerian. Here $ g_{\mid U_i} = \imath_{i}^{*}g $ where $\imath_{i} : U_{i}\longrightarrow M^{2n}$ is the inclusion. Also $(M^{2n},J,g)$ is globally conformal  K\"{a}ehler  if there is a $C^{\infty}$ function $f : M^{2n} \longrightarrow \mathbb{R}$ so that the metric $ exp(f)g $ is K\"{a}ehlerian \cite{Dr}.\\The first study
 on locally conformal K\"{a}ehler manifolds was done by Libermann in $1955$ \cite{L}. Visman \cite{v}, put down some geometrical conditions for locally conformal K\"{a}ehler manifold  and in $1982$ Tricerri mentioned different examples about the locally conformal K\"{a}ehler manifold \cite{Tr}. In $2001$, Banaru \cite{Ba} succeed to classify the sixteen classes of almost Hermitian manifold by using the two tensors of
Kirichenko, which are called the Kirichenko's tensors. Abood studied the properties of these tensors \cite{aa}. The locally conformal K\"{a}ehler manifold is one of the sixteen classes of almost Hermitian manifold. \\
In 1972, K. Kenmotsu introduced a class of contact metric manifolds, called Kenmotsu manifold, which is not Sasakian \cite{k}. The close relation between K\"{a}ehler manifolds and kenmotsu manifolds naturally leads to the question which objects, methods and theorems can be transfered from one to the other.\\
In this paper we get the idea of constructing globally conformal K\"{a}ehler manifolds and define conformal Kenmotsu manifolds which are not Kenmotsu, but by a conformal transformation the conformal Kenmotsu manifold is Kenmotsu. The paper is organized as follows:\\
In section 2, we introduce conformal Kenmotsu manifold and conformal Kenmotsu space form. In section 3, some definitions and results about conformal Kenmotsu space form and their CR-hypersurfaces are given. In section 4, 5, 6 and 7 We consider CR-hypersurfaces of a conformal Kenmotsu space form whose shape operator is parallel, scalar, recurrent or Lie $ \xi $-parallel, we show that if the Lee vector field of a conformal Kenmotsu space form is tangent  and normal to these type CR-hypersurfaces then the CR-hypersurfaces are totally geodesic and totally umbilic, respectively.

\section{riemannian geometry of conformal kenmotsu space form}
A 2n+1-dimensional differentiable manifold $ M $ is said to be an almost
contact metric manifold if it admits an almost contact metric structure $ (\varphi ,\tilde{\xi},\tilde{\eta},\tilde{g}) $ consisting of a tensor field $ \varphi $ of type (1,1), a vector field $\tilde{\xi}$, a 1-form
$\tilde{\eta}$ and a Riemannian metric $ \tilde{g} $ compatible with $ (\varphi ,\tilde{\xi},\tilde{\eta},\tilde{g}) $ satisfying
 following properties
\begin{eqnarray*}
 \varphi^{2}= -Id + \tilde{\eta} \otimes \tilde{\xi},\hspace{1cm}\tilde{\eta}(\tilde{\xi}) = 1,\hspace{1cm}\varphi\tilde{\xi} = 0,\hspace{1cm}\tilde{\eta} o\varphi =
 0,
 \end{eqnarray*}
  \begin{eqnarray*}
  \tilde{g}(\varphi X,\varphi Y)=\tilde{g}(X,Y) - \tilde{\eta}(X)\tilde{\eta}(Y),\hspace{1cm}\tilde{\eta}(X) =
  \tilde{g}(X,\tilde{\xi}),
 \end{eqnarray*}
 for all vector fields $ X, Y $ on $ M $.\\
An almost contact metric manifold $(M^{2n+1},\varphi ,\tilde{\xi},\tilde{\eta},\tilde{g}) $ is said to be Kenmotsu manifold [6] if
 \begin{eqnarray}
  &&(\tilde{\nabla}_{X}\varphi)Y= -\tilde{g}(X,\varphi Y)\tilde{\xi}-\tilde{\eta}(Y) \varphi X.
 \end{eqnarray}
 From $ (2.1) $, we have
 \begin{eqnarray}
 && \tilde{\nabla}_{X}\tilde{\xi}= X - \tilde{\eta}(X)\tilde{\xi},
 \end{eqnarray}
where, $\tilde{\nabla}$ denotes the Riemannian connection of $ \tilde{g} $. \\
Let $M$ be a Kenmotsu manifold. A plan section in $T_{m}M^{2n + 1}$ is called a $\varphi$-section if there exists a vector $X \in T_{m}M^{2n + 1}$ orthogonal to $ \xi $ such that $\{X,\varphi X\}$ span the section. The sectional curvature $K(X,\varphi X)$ is called $\varphi$-holomorphic sectional curvature. A Kenmotsu manifold $M$ with constant $\varphi$-holomorphic sectional curvature $c$ is said to be a Kenmotsu space form and is denoted by $M(c)$.
The necessary and sufficient condition for $M$ to have constant $\varphi$-holomorphic sectional curvature $c$ is
\begin{eqnarray}
\tilde{R}(X,Y)Z &=& \frac{c - 3}{4}\{\tilde{g}(Y,Z)X - \tilde{g}(X,Z)Y\}\nonumber\\
 &+& \frac{c + 1}{4} \{[\tilde{\eta}(X)Y - \tilde{\eta}(Y)X]\tilde{\eta}(Z)\nonumber\\
 & +& [\tilde{\eta}(Y)\tilde{g}(X,Z) - \tilde{\eta}(X)\tilde{g}(Y,Z)]\tilde{\xi}\nonumber\\
  &+& \tilde{g}(\varphi Y,Z)\varphi X - \tilde{g}(\varphi X,Z)\varphi Y - 2\tilde{g}(\varphi X,Y)\varphi Z\},
\end{eqnarray}
for all $X, Y, Z$ on $M$, where $ \tilde{R} $ is curvature tensor on $ M $ [6].\\
A  smooth manifold $M^{2n+1}$ with almost contact metric structure $(\varphi,\eta,\xi,g)$ is called a
conformal Kenmotsu manifold if there is a positive smooth function
$f:{M}^{2n+1}\rightarrow \mathbb{R}$ so that 
\begin{eqnarray*}
\tilde{g}=exp(f)g,\hspace{1cm}\tilde{\xi}={exp(-f)}^{\frac{1}{2}}\xi,
\hspace{1cm}\tilde{\eta}={exp(f)}^{\frac{1}{2}}\eta , \hspace{1cm}\tilde{\varphi} = \varphi
\end{eqnarray*}
is a Kenmotsu structure on $M^{2n+1}$.
Manifold $(M^{2n+1},\varphi,\eta,\xi,g)$ is called a conformal Kenmotsu space form if $M^{2n+1}$ with new almost contact metric structure $(\varphi,\tilde{\eta},\tilde{\xi},\tilde{g})$ is a Kenmotsu space form and is denoted by $ M(c) $ that $ c $ showes constant $\varphi$-holomorphic sectional curvature of Kenmotsu manifold  $(M^{2n+1},\varphi,\tilde{\eta},\tilde{\xi},\tilde{g})$.\\ Let $M^{2n+1}$ be a  conformal Kenmotsu manifold. Suppose $\tilde{\nabla}$ and $\nabla$
denote Riemannian connections $M^{2n+1}$ with respect to metrics
$\tilde{g}$ and $g$, respectively. Using Koszul formula, we
obtain the following relation between the connections
 $\tilde{\nabla}$ and $\nabla $ :
\begin{eqnarray}
\tilde{\nabla}_{X}Y = \nabla_{X}Y + \frac{1}{2}\{\omega(X)Y + \omega(Y)X - g(X,Y)\omega^{\sharp}\},
\end{eqnarray}
such that $\omega(X)= g(grad f,X) = X(f)$ and $g(\omega^{\sharp},X)=\omega(X)$ that $ \omega^{\sharp} = gradf $ is  the Lee vector field of $(M^{2n+1},\varphi,\eta,\xi,g)$.\\
Let $\tilde{R}$ and $R$ denote the curvature tensors on $(M^{2n+1},\varphi,\tilde{\eta},\tilde{\xi},\tilde{g})$ and $(M^{2n+1},\varphi,\eta,\xi,g)$, respectively.  Then the relation between $\tilde{R}$ and $R$ is given by
\begin{eqnarray}
\exp(-f)\tilde{g}(\tilde{R}(X,Y)Z,W)&=&g(R(X,Y)Z,W)\nonumber\\
&+&\frac{1}{2}\{B(X,Z)g(Y,W)-B(Y,Z)g(X,W)\nonumber\\
&+&B(Y,W)g(X,Z)-B(X,W)g(Y,Z)\}\nonumber\\
&+&\frac{1}{4}\|\omega^{\sharp}\|^{2}\{g(X,Z)g(Y,W)-g(Y,Z)g(X,W)\},
\end{eqnarray}
for all $ X, Y, Z $ on $ M $, where
\begin{eqnarray}
B:=\nabla \omega-\frac{1}{2}\omega\otimes\omega.
\end{eqnarray}
Furthermore, by the relations (2.1), (2.2) and (2.4) we get
\begin{eqnarray}
(\nabla_{X}\varphi)Y&=&(\exp(f))^{\frac{1}{2}}\{-{g(X,\varphi Y)\xi-\eta(Y)\varphi X}\}\nonumber\\
&-& \frac{1}{2}\{\omega(\varphi Y)X - \omega(Y)\varphi X+g(X,Y)\varphi \omega^{\sharp} - g(X,\varphi Y)\omega^{\sharp}\},\\
\nabla_{X}\xi &=&(\exp(f))^{\frac{1}{2}}\{X - \eta(X)\xi\}-\frac{1}{2}\{\omega(\xi)X - \eta(X)\omega^{\sharp}\}.
\end{eqnarray}
for all $ X, Y $ on $(M,\varphi,\eta,\xi,g)$.
\section{cr-hypersurfaces in conformal kenmotsu space form}
Let $M$ be $(2n+1)$-dimensional conformal Kenmotsu space form with almost contact metric structure $(\varphi,\xi,\eta,g)$. Consider a $2n$-dimensional manifold $\acute{M}$ embedded in $M$, $\acute{M}$ is called to be a $CR$-hypersurface of $M$ if the structure vector field $\xi$ is tangent to $\acute{M}$ and there exists
a pair of orthogonal complementary distribution $D$ and $D^{\bot}$ of $T\acute{M}$ such that \\
1) $D$ is invariant by $\varphi$, i.e. $\varphi(D_{p}) \subset D_{p}$, for any $p \in \acute{M}$;\\
2) $D^{\bot}$ is anti-invariant by $\varphi$, i.e. $\varphi(D^{\bot}_{p}) \subset T_{p}^{\bot}(\acute{M})$, for any $p \in \acute{M}$.\\
Now,  assumed $ \acute{M} $ be a hypersurface of a conformal Kenmotsu space form $M$ that for each $p \in \acute{M}$ the tangent vector $\xi$ always belong to the tangent bundle of the hypersurface $\acute{M}$. Let $\acute{g}$ is induced metric on $\acute{M}$, such that
\begin{eqnarray*}
\acute{g}(X,Y) = g(X,Y)
\end{eqnarray*}
for all $ X,Y \in T\acute{M} $. Also let $ N $ be the unit normal vector field to $\acute{M}$, we put $ \varphi N = -U $. Clearly $U$ is a unit tangent vector field on $\acute{M}$.
We denote by $D^{\bot} = span\{U, \xi\}$ the $2$-dimensional distribution generated by $U$ and $\xi$, and by $D$ the orthogonal complement of $D^{\bot}$ in $T\acute{M}$.
Thus we have the following decompositions
\begin{eqnarray}
 &&TM = D \oplus D^{\bot} \oplus span\{N\}\\
 &&T\acute{M} = D \oplus D^{\bot}.
\end{eqnarray}
Hence $ \acute{M} $ is a $ CR $-hypersurfase of $ M $.\\
Let $\acute{M}$ be a $CR$-hypersurface of a conformal Kenmotsu space form $M$. Denote by $\nabla$ and $\acute{\nabla}$ the Levi-Civita connection  on $M$ and induced Levi-Civita connection on $\acute{M}$ respectively, by using $(3.1)$ and $(3.2)$, the Gauss and Wiengarten formulas are
\begin{eqnarray*}
&&\nabla_{X}Y = \acute{\nabla}_{X}Y + h(X,Y) \\
&&\nabla_{X}N = -AX, \hspace{1cm}             \forall X,Y \in T\acute{M}
\end{eqnarray*}
where $A$ is the shape operator of $\acute{M}$ with respect to unit normal vector field $N$. It is known that
\begin{equation*}
h(X,Y) = \acute{g}(AX,Y)N,  \hspace{1cm}              \forall X,Y \in T\acute{M}.
\end{equation*}
In the usual way, by using $ (2.3) $ and $ (2.5) $ we derive the Codazzi equation;
\begin{eqnarray}
 && (\acute{\nabla}_{X}A)Y - (\acute{\nabla}_{Y}A)X = \frac{c + 1}{4}\exp(f) \{\acute{g}(X,U)\varphi Y - \acute{g}(Y,U)\varphi X \\
 &-& 2g(\varphi X,Y)U\} + \frac{1}{2}\{B(X,N)Y - B(Y,N)X\}\nonumber.
\end{eqnarray}
for all $X,Y \in T\acute{M}$.
\begin{lem}
Let $ \acute{M} $ be a $ CR $-hypersurface of a conformal Kenmotsu manifold $ (M,\varphi ,\xi,\eta,g)  $. Suppose that $ \omega^{\sharp} = N $. Then
\begin{eqnarray}
&&B(X,N) = 0,\\
&&B(X,Y) = -\acute{g}(AX,Y),\\
&&A\xi =\frac{1}{2}\xi ,\\
&&\acute{\nabla}_{X}\xi = {exp(f)}^{\frac{1}{2}}\{X - \eta(X)\xi \},
\end{eqnarray}
for all $ X, Y $ on $ \acute{M}$.\\
\begin{proof}
The lemma results from the Gauss and Weingarten formulas and relations (2.6) and (2.8).
\end{proof}
\end{lem}
\begin{lem}
Let $ \acute{M} $ be a $ CR $-hypersurface of a conformal Kenmotsu manifold $ (M,\varphi ,\xi,\eta,g)  $.  Let vector field $ \omega^{\sharp} $ is tangent to  $ \acute{M} $ and non-zero. Then
\begin{eqnarray}
&&B(X,N) = \omega(AX),\\
&&B(X,Y) = \acute{g}(\acute{\nabla}_{X}\omega^{\sharp},Y) - \frac{1}{2}\omega(X)\omega(Y),\\
&&A\xi = 0 ,\\
&&\acute{\nabla}_{X}\xi = {exp(f)}^{\frac{1}{2}}\{X - \eta(X)\xi \} -  \frac{1}{2}\{\omega(\xi)X - \eta(X)\omega^{\sharp} \},\hspace{1cm}
\end{eqnarray}
for all $ X, Y $ on $ \acute{M}$.
\begin{proof}
The lemma results from the Gauss and Weingarten formulas and relations (2.6) and (2.8).
\end{proof}
\end{lem}
\begin{lem}
Let $ \acute{M} $ be a $ CR $-hypersurface of a conformal Kenmotsu space form $ (M(c),\varphi ,\xi,\eta,g)  $. Suppose that $ \omega^{\sharp} = N $. Then
\begin{eqnarray}
&&\acute{R}(X,Y)\xi = exp(f)\{\eta(X)Y - \eta(Y)X \},
\end{eqnarray}
for all $ X, Y $ on $ \acute{M} $.
\begin{proof}
The lemma follows of the Gauss equation and relations (2.3), (2.5), (3.5) and (3.6).
\end{proof}
\end{lem}
\begin{lem}
Let $ \acute{M} $ be a $ CR $-hypersurface of a conformal Kenmotsu space form $ (M(c),\varphi ,\xi,\eta,g)  $. Let vector field $ \omega^{\sharp} $ is tangent to  $ \acute{M} $. Then
\begin{eqnarray}
&&\acute{R}(X,Y)\xi = (exp(f) - \frac{1}{4}\parallel \omega^{\sharp}\parallel^{2})\{\eta(X)Y - \eta(Y)X \}\\
&-&\frac{1}{2}\{\eta(\acute{\nabla}_{X}\omega^{\sharp})Y - \frac{1}{2}\omega(X)\omega(\xi)Y - \eta(\acute{\nabla}_{Y}\omega^{\sharp})X\nonumber\\
&+&\frac{1}{2}\omega(Y)\omega(\xi)X + \eta(X)\acute{\nabla}_{Y}\omega^{\sharp} - \frac{1}{2}\eta(X)\omega(Y)\omega^{\sharp}\nonumber\\
&-&   \eta(Y)\acute{\nabla}_{X}\omega^{\sharp} + \frac{1}{2}\eta(Y)\omega(X)\omega^{\sharp}\},\nonumber
\end{eqnarray}
for all $ X, Y $ on $ \acute{M} $.
\begin{proof}
The lemma follows of the Gauss equation and relations (2.3), (2.5), (3.9) and (3.10).
\end{proof}
\end{lem}
\section{cr-hypersurfaces  with parallel shape operator}
In [7], Kobayashi showed that a submanifold of a Kenmotsu manifold has parallel second fundamental form if and only if the submanifold is totally geodesic. As a generalization of this result we state the following results:
\begin{thm}\label{a}
 Let $ (M(c),\varphi ,\xi,\eta,g)$ is a conformal Kenmotsu space form. Suppose $ c\neq -1 $. Then there do not exist any CR-hypersurfaces $ \acute{M}$ in $ M$ with parallel shape operator.
\begin{proof}
Suppose that $\acute{\nabla}_{X}A = 0$ for all $X \in T\acute{M}$. If we put $ 0\neq X \in D $  and $ Y = U $ in $ (3.3) $, we have
\begin{equation}
-\frac{c + 1}{4}exp(f)\varphi X - \frac{1}{2}B(X,N)U + \frac{1}{2}B(U,N)X= 0.
\end{equation}
The set $ \{\varphi X, X, U\} $ is linearly independent. Then $ c = -1 $ which contradicts the hypothesis.
\end{proof}
\end{thm}
\begin{cor}
Let $\acute{M}$ be a CR-hypersurface of a conformal Kenmotsu space form $ (M(c),\varphi ,\xi,\eta,g)  $ with  shape operator $ A $. If $ A $ is parallel then
\begin{eqnarray*}
g(R(X,Y)Z,W) &=& (\frac{1}{4}\parallel \omega^{\sharp}\parallel^{2} - exp(f))\{g(Y,Z)X - g(X,Z)Y\} \\
&-& \frac{1}{2}\{B(X,Z)g(Y,W) -  B(Y,Z)g(X,W)\\
& +&  B(Y,W)g(X,Z) -  B(X,W)g(Y,Z)\},
\end{eqnarray*}
for all $ X, Y, Z $ and $ W $ on $ M $.
\begin{proof}
Since the shape operator $ A $ is parallel, by using theorem \ref{a} we get $ c = -1 $.  Now the above relation results from (2.3) and (2.5).
\end{proof}
\end{cor}
\begin{cor}
Let $\acute{M}$ be a CR-hypersurface of a conformal Kenmotsu space form $ (M(c),\varphi ,\xi,\eta,g)  $ with parallel shape operator. If vector field $ \omega^{\sharp} $ is tangent to $\acute{M}$ then $ A\omega^{\sharp} = 0 $.
\begin{proof}
Since the shape operator $ A $ is parallel, by using (4.1) we have
\begin{equation*}
- \frac{1}{2}B(X,N)U + \frac{1}{2}B(U,N)X= 0,
\end{equation*}
for all $ 0\neq X \in D $. The set $ \{X, U\} $ is linearly independent then $ B(X,N) = 0 $. Hence from (3.8) we get $ A\omega^{\sharp} = 0 $.
\end{proof}
\end{cor}
\begin{cor}
Let $\acute{M}$ be a CR-hypersurface of a conformal Kenmotsu space form $ (M(c),\varphi ,\xi,\eta,g)  $ with parallel shape operator. Suppose vector field $ \omega^{\sharp} $ is tangent to $\acute{M}$, then  $\acute{M}$ is   totally geodesic.
\begin{proof} 
Since the shape operator $ A $ is parallel we have
\begin{eqnarray}
&&\acute{R}(X,Y)AZ = \acute{\nabla}_{X}\acute{\nabla}_{Y}(AZ) - \acute{\nabla}_{Y}\acute{\nabla}_{X}(AZ) - \acute{\nabla}_{[X,Y]}(AZ) = A\acute{R}(X,Y)Z,
\end{eqnarray}
for all $ X, Y, Z $ on $\acute{M}$. Taking $ Y = Z = \xi $ in (4.2) and using from (3.10), (3.13) and corollary 4.3 we get
\begin{eqnarray}
&&(-exp(f) + \frac{1}{4}\parallel \omega^{\sharp}\parallel^{2} + \frac{1}{2}\eta(\acute{\nabla}_{\xi}\omega^{\sharp})X - \frac{1}{4}\omega(\xi)^{2})AX = 0,
\end{eqnarray}
for any $ X $ orthogonal to vector field $ \xi $.  Since $ \omega^{\sharp} $ is tangent to $\acute{M}$ we put $ \omega^{\sharp} = \alpha \zeta + \beta \xi$  where $ \zeta \in D\oplus \{U\}$ and $ \alpha,\beta $ are non-zero constants. Then from (3.11) and (4.3) we have $  exp(f)AX = 0$ for any $ X $ orthogonal to vector field $ \xi $. Hence  from (3.10) we get $ AX = 0 $ for any $ X $ on $\acute{M}$.
\end{proof}
\end{cor}
\begin{cor}
Let $\acute{M}$ be a CR-hypersurface of a conformal Kenmotsu space form $ (M(c),\varphi ,\xi,\eta,g)  $ with parallel shape operator. Suppose vector field $ \omega^{\sharp} = N $  then  $\acute{M}$ is totally umbilic.
\begin{proof}
Since the shape operator $ A $ is parallel we have (4.2). Taking $ Y = Z = \xi $ in (4.2) and using from (3.6) and (3.12) we obtain
\begin{eqnarray}
&&AX =  \frac{1}{2}X,
\end{eqnarray}
for any $ X $ orthogonal to vector field $ \xi $. (3.6) and (4.4) show $ AX = \frac{1}{2}X $ for any $ X $ on $\acute{M}$.
\end{proof}
\end{cor}
\section{cr-hypersurfaces  with scalar shape operator}
\begin{thm}\label{b}
Let $ (M(c),\varphi ,\xi,\eta,g)$ is a conformal Kenmotsu space form. Suppose $ c\neq -1 $. Then there do not exist any CR-hypersurfaces $ \acute{M}$ in $ M$ with scalar shape operator.
\begin{proof}
Suppose that $AX = \lambda X$ for all $X \in T\acute{M}$ where $ \lambda $ is a scalar function. By using  $ (3.3) $, we have
\begin{eqnarray*}
&&(X\lambda)Y - (Y\lambda)X = \frac{c + 1}{4}\exp(f) \{\acute{g}(X,U)\varphi Y - \acute{g}(Y,U)\varphi X \\
 &-& 2g(\varphi X,Y)U\} + \frac{1}{2}\{B(X,N)Y - B(Y,N)X\},
\end{eqnarray*}
for all $ X, Y $ on $\acute{M}$. If we put $ 0\neq X \in D $  and $ Y = U $ in above equation, we get
\begin{equation}
(X\lambda)U - (U\lambda)X = -\frac{c + 1}{4}exp(f)\varphi X - \frac{1}{2}B(X,N)U + \frac{1}{2}B(U,N)X.
\end{equation}
The set $ \{\varphi X, X, U\} $ is linearly independent. Then $ c = -1 $ which contradicts the hypothesis.
\end{proof}
\end{thm}
\begin{cor}
Let $\acute{M}$ be a CR-hypersurface of a conformal Kenmotsu space form $ (M(c),\varphi ,\xi,\eta,g)  $ with  shape operator $ A $. If $ A $ is scalar then
\begin{eqnarray*}
g(R(X,Y)Z,W) &=& (\frac{1}{4}\parallel \omega^{\sharp}\parallel^{2} - exp(f))\{g(Y,Z)X - g(X,Z)Y\} \\
&-& \frac{1}{2}\{B(X,Z)g(Y,W) -  B(Y,Z)g(X,W)\\
& +&  B(Y,W)g(X,Z) -  B(X,W)g(Y,Z)\},
\end{eqnarray*}
for all $ X, Y, Z $ and $ W $ on $ M $.
\begin{proof}
Since the shape operator $ A $ is scalar, by using theorem \ref{b} we get $ c = -1 $.  Now the above relation results from (2.3) and (2.5).
\end{proof}
\end{cor}
\begin{thm}
Let $\acute{M}$ be a CR-hypersurface of a conformal Kenmotsu space form $ (M(c),\varphi ,\xi,\eta,g)  $ with scalar shape operator. If vector field $ \omega^{\sharp} $ is tangent to $\acute{M}$ then $\acute{M}$ is totally geodesic.
\begin{proof}
It immediately results from (3.10).
\end{proof}
\end{thm}
\begin{thm}
Let $\acute{M}$ be a CR-hypersurface of a conformal Kenmotsu space form $ (M(c),\varphi ,\xi,\eta,g)  $ with scalar shape operator $ \lambda $. If vector field $ \omega^{\sharp} = N $,  then $\lambda = \frac{1}{2}$.
\begin{proof}
It immediately results from (3.6).
\end{proof}
\end{thm}
\section{cr-hypersurfaces  with recurrent shape operator}
\rm The notion of recurrent tensor field of type (r,s) on a differentiable manifold $ M $ with a linear connection was introduced in [8] and [13]. A non-zero tensor field $ K $ of type (r,s) on $ M $ which is said to be recurrent if there exists a 1-form $ \alpha $ such that $ \nabla K = \alpha \otimes K $. 
In [10], Sular and Ozgur showed a submanifold of a Kenmotsu manifold $ M $ has recurrent second fundamental form if and only if $ M $ is totally geodesic. As a generalization of this result we state the following results:\\
Let $\acute{M}$ be a CR-hypersurface of a conformal Kenmotsu space form $ (M(c),\varphi ,\xi,\eta,g)  $. $\acute{M}$ is said to be with recurrent shape operator if there exists a $ 1 $-form $ \alpha $ on $\acute{M}$ such that
\begin{eqnarray}
&&(\acute{\nabla}_{X}A)Y = \alpha(X)AY,
\end{eqnarray}
for all $ X, Y $ on $\acute{M}$. 
\begin{thm}
Let $\acute{M}$ be a CR-hypersurface of a conformal Kenmotsu space form $ (M,\varphi ,\xi,\eta,g)  $ with recurrent shape operator $ A $. Suppose vector field $ \omega^{\sharp} = N $.  Then $\acute{M}$ is totally umbilic.
\begin{proof}
If we put $ Y = \xi $ in (6.1) then from (3.6) and (3.7) it follows that
\begin{eqnarray*}
&&\frac{1}{2}exp(f)^{\frac{1}{2}}X - exp(f)^{\frac{1}{2}}AX = \frac{1}{2}\alpha(X)\xi ,
\end{eqnarray*}
for all $ X $ on $\acute{M}$. From inner product the above equation with vector field $ \xi $ and by using from (3.6) we obtain
\begin{eqnarray*}
&&\frac{1}{2}\alpha(X) =0,
\end{eqnarray*}
for any $ X $ on $\acute{M}$. Hence $ \alpha = 0 $. By using corollary 4.5 we get, $\acute{M}$ is totally umbilic.
\end{proof}
\end{thm}
\begin{thm}
Let $\acute{M}$ be a CR-hypersurface of a conformal Kenmotsu space form $ (M(c),\varphi ,\xi,\eta,g)  $ with recurrent shape operator $ A $. Suppose vector field $ \omega^{\sharp}  $ is tangent to $\acute{M}$.  Then either $\acute{M}$ is   totally geodesic or $ f = 2(Ln\beta - Ln2) $.
\begin{proof}
Let we put $ Y = \xi $ in (6.1), then from (3.10) and (3.11) we get
\begin{eqnarray}
&&(\frac{1}{2}\omega(\xi) - exp(f)^{\frac{1}{2}})AX = 0,
\end{eqnarray}
for any $ X $ orthogonal to $ \xi $. Since $ \omega^{\sharp} $ is tangent to $\acute{M}$ we put $ \omega^{\sharp} = \alpha \zeta + \beta \xi$  where $ \zeta \in D\oplus \{U\}$ and $ \alpha,\beta $ are non-zero constants. Then from (6.2) we have
\begin{eqnarray*}
&&(\frac{1}{2}\beta - exp(f)^{\frac{1}{2}})AX = 0,
\end{eqnarray*}
for any $ X $ orthogonal to $ \xi $. The above equation and (3.10) complete the proof of the lemma.
\end{proof}
\end{thm}
\section{cr-hypersurfaces  with Lie $ \xi $-parallel shape operator}
\rm Let $\acute{M}$ be a CR-hypersurface of a conformal Kenmotsu space form $ (M(c),\varphi ,\xi,\eta,g)  $. $\acute{M}$ is said to be with Lie $ \xi $-parallel shape operator if
\begin{eqnarray}
&&(L_{\xi}A)X = 0,
\end{eqnarray}
for any $ X $ on $\acute{M}$, where $ L_{\xi} $ shows Lie derivative along the direction $ \xi $.
\begin{thm}
Let $\acute{M}$ be a CR-hypersurface of a conformal Kenmotsu space form $ (M(c),\varphi ,\xi,\eta,g)  $ with Lie $ \xi $-parallel shape operator. Suppose vector field $ \omega^{\sharp} = N $.  Then $\acute{M}$ is totally umbilic.
\begin{proof}
 From (7.1) we have
\begin{eqnarray}
&&(\acute{\nabla}_{\xi}A)X - \acute{\nabla}_{AX}\xi + A\acute{\nabla}_{X}\xi = 0,
\end{eqnarray}
for any $ X $ on $\acute{M}$. From (3.3), (3.4) and (3.6) we obtain
 \begin{eqnarray}
&&(\acute{\nabla}_{\xi}A)X = \frac{1}{2}\acute{\nabla}_{X}\xi - A\acute{\nabla}_{X}\xi.
\end{eqnarray}
Now from (3.7), (7.2) and (7.3) we have
\begin{eqnarray}
&&exp(f)^{\frac{1}{2}}(\frac{1}{2}X - AX) = 0,
\end{eqnarray}
for any $ X $ on $\acute{M}$. (7.4) complete the proof of the lemma.
\end{proof}
\end{thm}
\begin{thm}
Let $\acute{M}$ be a CR-hypersurface of a conformal Kenmotsu space form $ (M,\varphi ,\xi,\eta,g)  $ with Lie $ \xi $-parallel shape operator $ A $. Suppose vector field $ \omega^{\sharp}  $ is tangent to $\acute{M}$. Then either $\acute{M}$ is   totally geodesic or $ f = 2(ln\beta - ln2) $.
\begin{proof}
From (3.3), (3.8) and (3.10) we obtain
 \begin{eqnarray}
&&(\acute{\nabla}_{\xi}A)X =  - A\acute{\nabla}_{X}\xi + \frac{1}{2}\omega(AX)\xi,
\end{eqnarray}
for any $ X $ on $\acute{M}$. Now from (3.11), (7.2) and (7.5) we have
\begin{eqnarray}
&&( \frac{1}{2}\omega(\xi)  - exp(f)^{\frac{1}{2}}) AX +  \frac{1}{2}\omega(AX)\xi = 0.
\end{eqnarray}
From inner product (7.6) with $ \xi $ and by using (3.10) we obtain $ \omega(AX) = 0 $ for any $ X $ on $\acute{M}$. Hence $ A\omega^{\sharp} = 0 $. Then from (7.6) we get
\begin{eqnarray*}
&&( \frac{1}{2}\omega(\xi)  - exp(f)^{\frac{1}{2}}) AX = 0,
\end{eqnarray*}
for any $ X $ on $\acute{M}$. Since $ \omega^{\sharp} $ is tangent to $\acute{M}$ we put $ \omega^{\sharp} = \alpha \zeta + \beta \xi$  where $ \zeta \in D\oplus \{U\}$ and $ \alpha,\beta $ are non-zero constants. Then from the above relation we have
\begin{eqnarray*}
&&( \frac{1}{2}\beta  - exp(f)^{\frac{1}{2}}) AX = 0,
\end{eqnarray*}
 for any $ X $ on $\acute{M}$. The above equation  complete the proof of the lemma.
\end{proof}
\end{thm}


\end{document}